\documentclass{article}
\usepackage[utf8]{inputenc}

\usepackage{graphicx,graphics}

\usepackage{rotating}
\usepackage[twoside]{geometry}
\usepackage{epsfig}
\usepackage[usenames,dvipsnames,svgnames,table]{xcolor}
\usepackage{cmap}

\usepackage{graphicx,graphics}

\usepackage{hyperref}
\hypersetup{
pdftitle={Strichartz&GFD},            
colorlinks=true,      
linkcolor=black,         
citecolor=blue}

\usepackage{amssymb,amsfonts,amstext,amsthm}
\usepackage{mathrsfs}
\usepackage[intlimits]{amsmath}
\newcommand{\dom}{{\mathrm{dom\,}}}

\newtheorem{lemma}{\bf Lemma}[section]
\newtheorem{theorem}{\bf Theorem}[section]

\newtheorem{remark}{Remark}[section]

\newtheorem{example}{Example}[section]

\newcounter{SE}

\title{On the Fourier asymptotics of absolutely continuous measures with power-law singularities}
\author{M. Aloisio\thanks{Corresponding author. Email: moacir@ufam.edu.br}, S. L. de Carvalho, C. R. de Oliveira and E. Souza}
\date{May 2022}

\begin{document}

\maketitle

\begin{abstract} We prove sharp estimates on the time-average behavior of the squared absolute value of the Fourier transform of some absolutely continuous measures that may have power-law singularities, in the sense that their Radon-Nikodym derivatives diverge with a power-law order. We also discuss an application to spectral measures of finite-rank perturbations of the discrete  Laplacian. 
\end{abstract}  

\ 


\

\noindent{\bf  AMS classification codes}: 28A80 (primary), 42A85 (secondary).    

\renewcommand{\thetable}{\Alph{table}}


\section{Introduction}\label{sectIntrod}

\subsection{Contextualization}

The study of the long time behavior of the Fourier transform of fractal (spectral) measures and of the modulus of continuity of the distribution of such measures play an important role in spectral theory and quantum dynamics (such behavior is related to good transport properties). Actually, most of these works are motivated by possible applications to Schr\"odinger operators (see \cite{AvilaUaH,Last,strichartz1990,Zhao,Zhao2} and references therein). In this context, one may highlight two classical results on finite Borel measures on~$\mathbb{R}$: the Riemann-Lebesgue's Lemma (around 1900) and the Wiener's Lemma (around 1935).

One may also highlight Strichartz's Theorem \cite{strichartz1990}, from 1990 (see Theorem \ref{Strichartztheorem} (i) below), that establishes (power-law) convergence rates for the time-average behavior of the squared absolute value of the Fourier transform of uniformly $\alpha$-H\"older continuous measures. We present some details.

Let $\mu$ be a finite positive Borel measure on~$\mathbb{R}$ and $\alpha \in [0,1]$. We recall that $\mu$ is uniformly $\alpha$-H\"older continuous (denoted {\rm U}$\alpha${\rm H}) if there exists a constant $C>0$ such that for each interval $I$ with $\ell(I) < 1$, $\mu(I) \le C\, \ell(I)^\alpha$, where $\ell(\cdot)$ denotes the Lebesgue measure on~$\mathbb{R}$. Namely, one has the following result.

\begin{theorem}[Theorems 2.5 and 3.1  in \cite{Last}]\label{Strichartztheorem} Let $\mu$ be a finite Borel measure on $\mathbb{R}$ and $\alpha \in [0,1]$.

\begin{enumerate}

\item[\rm{i)}] If $\mu$ is {\rm U}$\alpha${\rm H}, then there exists a constant $C_\mu> 0$, depending only on $\mu$, such that for every $f \in {\mathrm L}^2(\mathbb{R}, d\mu )$ and every $t>0$, 
\[\frac{1}{t} \int_0^t \bigg|\int_{\mathbb{R}} e^{-2\pi isx} f(x)\,  d\mu(x) \bigg|^2  ds < C_{\mu} \|f\|_{{\mathrm L}^2(\mathbb{R}, d\mu )}^2 t^{-\alpha}. \] 

\item[\rm{ii)}] If there exists $C_{\mu}>0$ such that for every $t>0$, 
\[\frac{1}{t} \int_0^t \bigg|\int_{\mathbb{R}} e^{-2\pi isx}\,  d\mu(x) \bigg|^2  ds < C_{\mu}  t^{-\alpha},\] 
then $\mu$ is {\rm U}$\frac{\alpha}{2}${\rm H}.

\end{enumerate}
\end{theorem} 

\begin{remark}{\rm Note that Theorem \ref{Strichartztheorem}-i) is, indeed, a particular case of Strichartz's Theorem~\cite{strichartz1990}, which holds for $\sigma$-finite measures.}
\end{remark}

Motivated by applications in spectral theory and quantum dynamics, we use in this work Fourier analysis to prove sharp estimates on the time-average behavior of the squared absolute value of the Fourier transform of some absolutely continuous measures that may have power-law singularities. Our main goal here is to obtain initial states (for the Schr\"odinger equation) for which the respective spectral measures have a singularity with a power-law growth rate, and for which the asymptotic behavior of the respective (time-average) quantum return probabilities (see definition ahead) depends continuously on such singularities (see Theorem \ref{maintheorem} and Example \ref{ex2} ahead). To the best knowledge of the present authors, this phenomenon has never been discussed, although it may be natural to specialists.

In the next remark we discuss the fact that Theorem \ref{Strichartztheorem}-i), in general, is not sufficient to obtain  sharp estimates on such  Fourier transform averages of absolutely continuous measures with power-law singularities. 

\begin{remark}{\rm For some important classes of measures in spectral theory and quantum dynamics, such as spectral measures of dynamically defined Schr{\"o}dinger operators, $1/2$-H\"older continuity is typically optimal (usually due to the fact that there are square root singularities associated with the  boundary  of the spectrum; see \cite{AvilaUaH,Damanik,Zhao,Zhao2} for additional comments); by  Theorem~\ref{Strichartztheorem}-i), the time-average behavior of the squared absolute value of the Fourier transform of such measures decays  at least as $1/\sqrt{t}$. This rate, in general, is far from optimal as, for instance, is the case of the  discrete Laplacian:  let $\ell^2({\mathbb{Z}})$ and $\delta_j = (\delta_{jk})_{k \in {\mathbb{Z}}}$, $j\in{\mathbb Z}$, be its canonical basis, and consider the Laplacian, whose action on $\psi\in\ell^2({\mathbb{Z}})$ is given by  
\[\triangle\psi(k) = \psi(k+1) + \psi(k-1);\]
so, although the spectral measure $\mu_{\delta_0}^{\triangle}$ of the pair $(\triangle,\delta_0)$ is at most uniformly $1/2$-H\"older continuous,  one has
\begin{equation}\label{eq00}
\frac{1}{t} \int_0^t \bigg|\int_{\mathbb{R}} e^{-2\pi isx} d\mu_{\delta_0}^{\triangle}(x)\bigg|^2  ds=O(\log(t)/t);
\end{equation}
see the case $\beta = \frac{1}{2}$ in Example \ref{ex1} and, e.g., Section 12.3 in \cite{Oliveira} for details of the Radon-Nikodym derivative of such spectral measure (here, $h(t) = O(r(t))$ indicates that there is $C>0$ so that, for each $t>0$, $h(t)\leq Cr(t)$).  

By Theorem~\ref{Strichartztheorem}-i), one may conclude that
\begin{equation*}
\frac{1}{t} \int_0^t \bigg|\int_{\mathbb{R}} e^{-2\pi isx} d\mu_{\delta_0}^{\triangle}(x)\bigg|^2  ds=O(1/\sqrt{t}),
\end{equation*}
which gives a worse bound than the one given by~\eqref{eq00}. Moreover, since $\mu_{\delta_0}^{\triangle}$  is at most uniformly $1/2$-H\"older continuous, by Theorem \ref{Strichartztheorem}-ii), the rate in~(\ref{eq00}) is (power-law) optimal. Namely, suppose that there exists $\varepsilon>0$ such that one can replace $t^{-1}$ by $t^{-1-\varepsilon}$  in~\eqref{eq00}; then, by Theorem \ref{Strichartztheorem}-ii),  $\mu_{\delta_0}^{\triangle}$ is at least uniformly $(1/2+\varepsilon/4)$-H\"older continuous.}
\end{remark}

In order to put our work into perspective, we present the following example. To each $0<\beta <1$, denote by  
\begin{equation}\label{eqMbeta}
M_\beta:= \max_{\eta > 0}\bigg|\displaystyle\int_0^\eta e^{-iu} u^{-\beta} du \bigg|^2;
\end{equation}
for $0<\beta <1$,  $\bigg|\displaystyle\int_0^\infty e^{-iu} u^{-\beta} du \bigg| = \Gamma(1-\beta)$, where $\Gamma$ stands for the Gamma Function; thus, $M_\beta < \infty$.

\begin{example}\label{ex1} {\rm Set, for each $\frac{1}{2} \leq \beta < 1$ and each $f \in {\mathrm L}^1(\mathbb{R})$, 
\begin{equation}\label{eqKB}
K_{\beta,f} := \{(|\cdot|^{-\beta}  \chi_{(0,1]}) \ast f\},
\end{equation}
so $K_{\beta,f}  \in {\mathrm L}^1(\mathbb{R})$. By the Convolution Theorem, it follows that for each $s >0$,  
\begin{eqnarray*} \widehat{K_{\beta,f}}(s) &=& \biggl\{\int_0^1 e^{-2\pi i x s} x^{-\beta} dx \biggl\} \, \hat{f}(s) = \biggl\{\frac{1}{(2\pi)^{1-\beta}s^{1-\beta}} \int_0^1 e^{-2\pi i x s} (2 \pi x s)^{-\beta} (2\pi s)\,dx \biggl\} \, \hat{f}(s)\\ &=& \biggl\{\frac{1}{(2\pi)^{1-\beta}s^{1-\beta}} \int_0^{2\pi s} e^{- iu} u^{-\beta} du \biggl\} \, \hat{f}(s),
\end{eqnarray*}
and so, for each $t >0$ and $\frac{1}{2} < \beta < 1$,
\begin{eqnarray}\label{eqex}
  \nonumber\frac{1}{t} \int_0^t \bigg|\int_{\mathbb{R}} e^{-2\pi isx} d\mu_{\beta,f}(x)\bigg|^2  ds  &\leq& {M_\beta}\, \frac{1}{t}\int_0^ts^{2(\beta-1)}|\hat{f}(s)|^2 ds\\
    &\le& {M_\beta} \frac{\|f\|^2_{{\mathrm L^1}(\mathbb{R})}}{t}\int_0^ts^{2(\beta-1)}ds
  = \frac{{M_\beta}}{2\beta-1} \, \frac{\|f\|^2_{{\mathrm L^1}(\mathbb{R})}}{t^{2(1-\beta)}},
\end{eqnarray}
where $d\mu_{\beta,f}(x) =  K_{\beta,f}(x)\, dx$. For $\beta = \frac{1}{2}$ and for every $t> 1,$ one has
\begin{eqnarray}\label{eqex777777}
  \nonumber\frac{1}{t} \int_0^t \bigg|\int_{\mathbb{R}} e^{-2\pi isx} d\mu_{\beta,f}(x)\bigg|^2  ds &=&   \nonumber\frac{1}{t} \int_0^1 \bigg|\int_{\mathbb{R}} e^{-2\pi isx} d\mu_{\beta,f}(x)\bigg|^2  ds +   \nonumber \frac{1}{t} \int_1^t \bigg|\int_{\mathbb{R}} e^{-2\pi isx} d\mu_{\beta,f}(x)\bigg|^2  ds\\ \nonumber &\leq&  \frac{\|K_{\beta,f}\|_{{\mathrm L}^1(\mathbb{R})}^2}{t}  + \, \frac{M_\beta\|f\|^2_{{\mathrm L^1}(\mathbb{R})}}{t}\int_1^t \frac{1}{s} \, ds\\ &=&  \frac{\|K_{\beta,f}\|_{{\mathrm L}^1(\mathbb{R})}^2}{t}  + \, \frac{M_\beta\|f\|^2_{{\mathrm L^1}(\mathbb{R})} \log(t)}{t}.
\end{eqnarray}

We argue that the above  power-law upper estimates cannot be improved; suppose, on the contrary, that there exists $0<\varepsilon<1-\beta$ such that one can replace $2(1-\beta)$ by $2(1-\beta)+\varepsilon$ in the estimate~\eqref{eqex} or \eqref{eqex777777} (recall that in \eqref{eqex777777} $\beta = 1/2$). Then, by Theorem \ref{Strichartztheorem}-ii),  for each $f \in {\mathrm L}^1(\mathbb{R})$, $\mu_{\beta,f}$ is at least {\rm U}$(1-\beta+ \varepsilon/4)${\rm H}.

Now, set, for each $0<\delta <1-\beta$ and each $x \in \mathbb{R}\setminus\{0\}$, 
\[f_\delta(x) = \frac{1}{x^{1-\delta}} \; \chi_{(0,1]}(x),\]
and $f_\delta \in {\mathrm L}^1(\mathbb{R})$. Note that for  $0< x \leq 1$,
\begin{eqnarray*}
K_{\beta,f_\delta}(x) &=&  \int_0^1 \frac{1}{y^{\beta} |x-y|^{1-\delta}} \, dy \geq  \int_0^x \frac{1}{y^{\beta} |x-y|^{1-\delta}} \, dy \geq \frac{1}{x^\beta} \int_0^x \frac{1}{|x-y|^{1-\delta}} \, dy \geq    x^{- \beta+\delta}.
\end{eqnarray*}
Thus, for  $0< \epsilon < 1$,
\[ \int_0^\epsilon K_{\beta,f_\delta}(x) \, dx \geq \int_0^\epsilon x^{-\beta+\delta} dx  =  \frac{\epsilon^{(1- \beta+ \delta)} }{(1-\beta + \delta)},\]
and therefore, $\mu_{\beta,f_\delta}$ is at most  {\rm U}$(1-\beta+ \delta)${\rm H}.

Finally, let $0<\delta<\varepsilon/4$; then, since $\mu_{\beta,f_\delta}$ is at most {\rm U}$(1-\beta+ \delta)${\rm H}, one gets a contradiction with the fact that $\mu_{\beta,f_\delta}$ is also {\rm U}$(1-\beta+ \varepsilon/4)${\rm H}. 

We emphasize that by applying Theorem~\ref{Strichartztheorem}-i),  it may be obtained at most that 
\[\frac{1}{t} \int_0^t \bigg|\int_{\mathbb{R}} e^{-2\pi isx} d\mu_{\beta,f}(x)\bigg|^2  ds=O(t^{-(1-\beta)+\delta}),
\]
which gives a worse bound than the one given by~\eqref{eqex} (respec.~\eqref{eqex777777} for $\beta = \frac{1}{2}$) for $0<\delta<1-\beta$.

We also note that what makes this example interesting is the fact that the measure has a power-law singularity, in the sense that its Radon-Nikodym derivative has a power-law divergence. }
\end{example}

In this work we use Fourier analysis to extend the estimates in (\ref{eqex})-\eqref{eqex777777} to  the measures  
\[
d\mu_{\beta,f,g}(x) =  K_{\beta,f}(x) g(x)\, dx,
\]
with $f \in {\mathrm L}^1(\mathbb{R})$ and $g \in {\mathrm L}^\infty[0,1]$ (see Theorem~\ref{2Stheorem} ahead),  and then we discuss an application of this result to spectral measures of finite-rank perturbations of the  Laplacian. Namely, by taking into account Example~3.1 in~\cite{Last}, we use this class of measures to obtain initial states (for the Schr\"odinger equation) for which the respective spectral measures have power-law singularities, and for which the asymptotic behavior of the respective (time-average) quantum return probabilities (see definition ahead) depend continuously on such singularities (see Theorem~\ref{maintheorem} and Example~\ref{ex2}). 

The organization of this work is as follows. In Subsection~\ref{subsectStrich} we discuss a Strichartz's-like Inequality (Theorem~\ref{2Stheorem}). In Section~\ref{sectFRPFL} we use some well-known results on the Radon-Nikodym derivative of spectral measures \cite{Germinet,Lasttransfermatrix} to present an application to finite-rank perturbations of the  Laplacian.  The proof of Theorem~\ref{2Stheorem} is left to Section~\ref{sectProofMain}.

Some words about the notation:~$\hat{f}$  will always denote the Fourier transform of a function $f \in {\mathrm L}^1(\mathbb{R})$. If $h,g : \mathbb{R} \longrightarrow \mathbb{R}$ are mensurable functions, then $h \ast g$ denotes the convolution product of $h$ and $g$; $\mu$ always indicates  a finite positive Borel measure on $\mathbb{R}$. For each $x \in \mathbb{R}$ and each $\epsilon>0$,  $B(x,\epsilon)$ denotes the open interval $(x-\epsilon,x+\epsilon)$. If $g$ is a complex-valued function, then $\mathfrak{Re}(g)$ and $\mathfrak{Im}(g)$ denote its real and the imaginary parts, respectively. If~$f$ is a real-valued function, then $f^+$ and $f^-$ denote its positive and the negative parts, respectively.

\subsection{A Strichartz's-like Inequality} \label{subsectStrich}

Let $K_{\beta,f}$ be as in Example \ref{ex1}, $g \in {\mathrm L}^\infty[0,1]$ and consider
\begin{equation}
d\mu_{\beta,f,g}(x) =   K_{\beta,f}(x) g(x) \, dx.
\end{equation}
For simplicity, suppose that $f,g$ are nonnegative (measurable) real-valued functions. So, by well-known arguments \cite{Last,strichartz1990} (see (\ref{maineq1}) ahead), it is possible to show that for every $t>0$, 
\begin{eqnarray} \label{eq0}
 \nonumber \frac{1}{t} \int_0^t \bigg|\int_{\mathbb{R}} e^{-2\pi isx} d\mu_{\beta,f,g}(x) \bigg|^2\nonumber   ds &=&  \nonumber \frac{1}{t} \int_0^t \bigg|\int_{\mathbb{R}} e^{-2\pi isx} K_{\beta,f}(x) g(x) \, dx \bigg|^2  ds\\ &\le&  \frac{e^{2\pi} }{2 \sqrt{\pi}}   \int_{\mathbb{R}} \int_{\mathbb{R}} K_{\beta,f}(x)g(x)  K_{\beta,f}(y) g(y)   e^{-\frac{t^2|x-y|^2}{4}} dx dy. 
\end{eqnarray}
Moreover, for each $x \in \mathbb{R}$ and each $0<\epsilon<1$, one has 
\[ \mu_{\beta,f}(B(x,\epsilon)) \leq \|f\|_{{\mathrm L}^1(\mathbb{R})} \epsilon^{1-\beta},\] 
where $d\mu_{\beta,f}(x) =  K_{\beta,f}(x) dx$. So, by using (\ref{eq0}) and a Strichartz's-like  argument (as in \cite{Last,strichartz1990}), it follows that for every $t>0$,
\begin{eqnarray} \nonumber \frac{1}{t} \int_0^t \bigg|\int_{\mathbb{R}} e^{-2\pi isx} d\mu_{\beta,f,g}(x) \bigg|^2   ds &\leq&    \frac{e^{2\pi} \|g\|_{{\mathrm L}^\infty[0,1]}^2}{2 \sqrt{\pi}}  \int_{\mathbb{R}} \int_{\mathbb{R}} K_{\beta,f}(x)  K_{\beta,f}(y)   e^{-\frac{t^2|x-y|^2}{4}} dx dy\\ \nonumber &=&  \frac{e^{2\pi} \|g\|_{{\mathrm L}^\infty[0,1]}^2}{2 \sqrt{\pi}}   \int_{\mathbb{R}} \int_{\mathbb{R}}   e^{-\frac{t^2|x-y|^2}{4}} d\mu_{\beta,f}(y)  d\mu_{\beta,f}(x)\\  &=&\nonumber     \frac{e^{2\pi} \|g\|_{{\mathrm L}^\infty[0,1]}^2}{ \sqrt{\pi}}   \int_{\mathbb{R}} \sum_{n=0}^\infty \int_{\frac{n}{t}\leq |x-y|<\frac{n+1}{t}}   e^{-\frac{t^2|x-y|^2}{4}} d\mu_{\beta,f}(y)  d\mu_{\beta,f}(x) \\ \nonumber  &\leq&    \frac{e^{2\pi} \|g\|_{{\mathrm L}^\infty[0,1]}^2}{ \sqrt{\pi}}   \int_{\mathbb{R}} \sum_{n=0}^\infty e^{-n^2/4} \|f\|_{{\mathrm L}^1(\mathbb{R})} t^{\beta-1}  d\mu_{\beta,f}(x)\\   &=&   \frac{e^{2\pi} \|g\|_{{\mathrm L}^\infty[0,1]}^2}{ \sqrt{\pi}}  \left(\sum_{n=0}^\infty e^{-n^2/4}\right) \|f\|_{{\mathrm L}^1(\mathbb{R})} \|K_{\beta,f}\|_{{\mathrm L}^1(\mathbb{R})}\; t^{-(1-\beta)}.
\end{eqnarray}
Thus, by Young's Convolution Inequality, one gets, for every $t>0$,
\begin{eqnarray}\label{eq2} \nonumber \frac{1}{t} \int_0^t \bigg|\int_{\mathbb{R}} e^{-2\pi isx} d\mu_{\beta,f,g}(x) \bigg|^2   ds  &\leq&   \frac{e^{2\pi} \|g\|_{{\mathrm L}^\infty[0,1]}^2}{ \sqrt{\pi}}     \left(\sum_{n=0}^\infty e^{-n^2/4}\right) \|f\|^2_{{\mathrm L}^1(\mathbb{R})} \|(|\cdot|^{-\beta} \chi_{(0,1]})\|_{{\mathrm L}^1(\mathbb{R})} t^{-(1-\beta)} \\ &=& \left(\sum_{n=0}^\infty e^{-n^2/4}\right) \frac{ e^{2\pi}  \|g\|_{{\mathrm L}^\infty[0,1]}^2 \|f\|_{{\mathrm L}^1(\mathbb{R})}^2 }{ \sqrt{\pi}(1-\beta)}\;    t^{-(1-\beta)}.  
\end{eqnarray}

We remark that the uniform estimate over $x$ in the discussion above makes the decay in~\eqref{eq2} far from optimal (see the proof of Theorem \ref{2Stheorem} and compare~\eqref{eq2} with~(\ref{eq00lemma}) and~(\ref{eq1lemma})). By using Fourier analysis, we will explore this point of the argument to obtain the following result.

\begin{theorem}\label{2Stheorem} For $\frac{1}{2} \leq \beta < 1$, let $ K_{\beta,f}$ and~$M_\beta$ be as before.  To every  $g \in {\mathrm L}^\infty[0,1]$, consider $d\mu_{\beta,f,g}(x)= K_{\beta,f}(x) g(x)  dx$. Then:

\begin{enumerate}
\item[\rm {i)}] if $\frac{1}{2} < \beta < 1$,  for every $t>0$, 
\begin{eqnarray*}
 \frac{1}{t} \int_0^t \bigg|\int_{\mathbb{R}} e^{-2\pi isx}\, d\mu_{\beta,f,g}(x) \bigg|^2  ds  &\leq& 2^{18} e^{2\pi} {M_\beta} \Gamma(\beta - 1/2)  \|f\|_{{\mathrm L}^1(\mathbb{R})}^2  \|g\|_{{\mathrm L}^\infty[0,1]}^2\;   t^{-2(1-\beta)};
\end{eqnarray*}

\item[\rm {ii)}] if $\beta = \frac{1}{2}$,  for every $t>0$, 
\begin{eqnarray*}
 \frac{1}{t} \int_0^t \bigg|\int_{\mathbb{R}} e^{-2\pi isx}\, d\mu_{\beta,f,g}(x) \bigg|^2  ds  &\leq& e^{2\pi} \|f\|_{{\mathrm L}^1(\mathbb{R})}^2  \|g\|_{{\mathrm L}^\infty[0,1]}^2 \biggr[\biggr( \frac{1}{t} +  M_{\frac{1}{2}}  \frac{\Gamma(0, 4\pi^2/t^2)}{t}\biggr) \biggr],
\end{eqnarray*}
in particular, for sufficiently large $t$ ,
\begin{eqnarray*}
 \frac{1}{t} \int_0^t \bigg|\int_{\mathbb{R}} e^{-2\pi isx}\, d\mu_{\beta,f,g}(x) \bigg|^2  ds  &\leq& e^{2\pi} \|f\|_{{\mathrm L}^1(\mathbb{R})}^2  \|g\|_{{\mathrm L}^\infty[0,1]}^2 \biggr[\biggr( \frac{1}{t} +  3 M_{\frac{1}{2}}  \frac{\log(t)}{t}\biggr) \biggr],
\end{eqnarray*}
 since
 \[\lim_{t \to \infty} \frac{\Gamma(0, 4\pi^2/t^2)}{\log(t)} = 2,\]
\end{enumerate}
where $\Gamma(\cdot,\cdot)$ denotes the Incomplete Gamma Function.
\end{theorem}
\begin{remark}
\end{remark}

\begin{enumerate} 

\item [i)] As mentioned in Example~\ref{ex1}, in general, one cannot get a better power-law estimate than $O(t^{-2(1-\beta)})$ for all $f \in {\mathrm L}^1(\mathbb{R})$.

\item [ii)]  If $0 \leq \beta < \frac{1}{2}$ then,  by Young's Convolution Inequality, $K_{\beta,f} \cdot g  \in  {\mathrm L}^2(\mathbb{R})$. Hence, by applying Theorem \ref{Strichartztheorem}-i) to $K_{\beta,f} \cdot g$ and to $\chi_{[0,1]} \, dx$ (which is  U$1$H), one gets
\begin{equation*}
 \frac{1}{t} \int_0^t \bigg|\int_{\mathbb{R}} e^{-2\pi isx} d\mu_{\beta,f,g}(x) \bigg|^2  ds = \frac{1}{t} \int_0^t \bigg|\int_{\mathbb{R}} e^{-2\pi isx} K_{\beta,f}(x)  \, g(x)dx \bigg|^2  ds = O(t^{-1}),
\end{equation*} and $\beta=1/2$ is a transition point, so justifying its peculiar behavior as in Theorem~\ref{2Stheorem}-ii).

\item [iii)]  For an arbitrary $h \in {\mathrm L}^1(\mathbb{R})$, it is well known that $\hat{h}(s)$ can decay arbitrarily slow (see, e.g.,~\cite{Muller}). In this context, for $\frac{1}{2} \leq \beta < 1$, it is particularly interesting that the power-law asymptotic behavior of 
\begin{eqnarray*}
 \frac{1}{t} \int_0^t \bigg|\int_{\mathbb{R}} e^{-2\pi isx} K_{\beta,f}(x)  \, g(x)dx \bigg|^2  ds
\end{eqnarray*}
is inherited from the asymptotic behavior of the Fourier transform of $|\cdot|^{-\beta}$, which depends continuously on $\beta$.
\end{enumerate}

\section{Finite-rank perturbations of the  Laplacian}\label{sectFRPFL}

Let $\delta_j = (\delta_{jk})_{k \in \mathbb{N}}$, $j=1,2...$, be the canonical basis of~$\ell^2(\mathbb{N})$. Consider the  Laplacian with  Dirichlet boundary condition, whose action on $\psi\in\ell^2(\mathbb{N})$ is 
\[\triangle\psi(k) = \psi(k+1) + \psi(k-1),\]
with $\psi(0)=0$. Consider also the finite-rank perturbations of the Laplacian with Dirichlet boundary condition, acting in $\ell^2(\mathbb{N})$ by the law
\[
\begin{cases} H_0 \psi \, \, =-\triangle\psi, \\ H_N \psi = -\triangle\psi  + \displaystyle\sum_{j=1}^N v_j \langle \psi, \delta_j \rangle \delta_j,\qquad N \geq 1,  \end{cases}
\]
with $\psi(0)=0$, where $(v_n)_{n \in {\mathbb{N}}}$ is a given real sequence.

The study of the dynamical and spectral properties of these operators is a classical subject in spectral theory for at least two reasons: it is a relatively simple model, and therefore it can be used to discuss results in quantum dynamics by avoiding technical complications; they can be used as approximations to some Schr\"odinger operators (see, for instance, Section~3 in~\cite{Germinet}). In this context, our main goal here is to study the quantum dynamics of these operators for some initial states whose spectral measures may have power-law singularities. Naturally, some preparation is required.  

Let $\mu_\psi^{N}$ denote the spectral measure associated with the pair $(H_N,\psi)$ and $\displaystyle\frac{d\mu_\psi^{N}}{dx}$  the Radon-Nikodym derivative of  $\mu_\psi^{N}$ with respect to the Lebesgue measure, and let $R_N(z) = (H_N - z)^{-1}$ be the corresponding resolvent operator. It is well known that for every $N \in \mathbb{N} \cup \{0\}$, $H_N$ has  purely absolutely continuous spectrum (for details, see \cite{Lasttransfermatrix}). Thus, in this case, $\displaystyle\frac{d\mu_\psi^{N}}{dx}\in {\mathrm L}^1(\mathbb{R})$. Recall that for every $N \in \mathbb{N} \cup \{0\}$, the transfer matrix $T_N(E,n,n-1)$ between the sites $n-1$ and $n$, $n\in\mathbb{N}$, is given by 
\[T_N(E,n,n-1) =\left(\begin{array}{cc}
    E-v_n\chi_{[1,N]}(n) & -1\\ 
   1 & 0\\
\end{array}\right).\]
Moreover, if $u_\theta(E,n)$, with $ \theta \in [0,\pi]$, denotes the solution to the eigenvalue equation $H_Nu = Eu$ at $E \in \mathbb{R}$ that satisfies 
\[u_\theta(E,m) = \sin(\theta), \, \, \, \, u_\theta(E,m+1) = \cos(\theta),\]
then 
\[\left(\begin{array}{cc}
    u_\theta(E,n+1)\\ 
   u_\theta(E,n)\\
\end{array}\right)= T_N(E,n,m)\left(\begin{array}{cc}
    u_\theta(E,m+1)\\ 
   u_\theta(E,m)\\
\end{array}\right),\]
with $T_N(E,n,m)=T_N(E,n,n-1)\cdots T_N(E,m+1,m)$. We need the following technical result.

\begin{lemma}\label{teclemma}  Let $N \in \mathbb{N} \cup \{0\}$. Then, there exist constants $C_{1,N},C_{2,N}> 0$ such that, for every $E \in [0,1]$, 
\[C_{1,N}\, \leq\, \displaystyle\frac{d\mu_{\delta_1}^{N}}{dx}(E)\, \leq\, C_{2,N};\]
in particular, $\displaystyle\frac{d\mu_{\delta_1}^{N}}{dx} \in {\mathrm L}^\infty[0,1]$. 
\end{lemma}

\begin{remark}{\rm The result stated in Lemma \ref{teclemma} is expected, since the boundary of the interval $[0,1]$ is far from $\pm 2$, which are the only points where $\displaystyle\frac{d\mu_{\delta_1}^{H_0}}{dx}$ diverges; see \cite{Oliveira} for details on the spectral measure $\mu_{\delta_1}^{H_0}$. Although natural to specialists, we present a proof of this result for the convenience of the reader.}
\end{remark}

\begin{proof} [{Proof} {\rm (Lemma~\ref{teclemma})}] Let $N\in\mathbb{N}\cup\{0\}$. It follows from Lemma 3.1 in \cite{Germinet} that there exists $D> 0$ such that for every $E \in [0,1]$,
  \[\displaystyle\frac{d\mu_{\delta_1}^{N}}{dx}(E) \geq  \frac{D}{\|T_N(E,N,0)\|^2}.\] 
Note that there exists $F_N>0$ such that for each $E\in[0,1]$, $\Vert T_N(E,N,0)\Vert<F_N$ (since $T_N(E,N,0)$ is the product of $N$ matrices whose norms are bounded); thus,

\[0 < C_{1,N} := \frac{D}{F_N^2} \leq \inf_{E \in [0,1]}\displaystyle\frac{d\mu_{\delta_1}^{N}}{dx}(E).\]

Now, if $n \geq N$,  then $T_N(E,n+1,n) =  T_0(E,n+1,n)=A(E)$, where
  \begin{eqnarray*} A(E)=\left(\begin{array}{cc} E&-1\\1&0\end{array}\right).
  \end{eqnarray*}
  It is straightforward to show that for $E\in[0,1]$, $A(E)$ is similar to a rotation matrix; thus, there exists $C_N> 0$ such that for each $E \in [0,1]$ and each $n \geq N$, $\|T_N(E,n,0)\|^2 \leq C_N$. Indeed,
  \[\|T_N(E,n,0)\|=\Vert T_N(E,n,N)\cdot \cdot \cdot T_N(E,N,0)\Vert\le \Vert A(E)^{n-N}\Vert\cdot\Vert T_N(E,N,0)\Vert;\]
since $A(E)^{n-N}$ is similar to a rotation and since for each $E\in[0,1]$, $\Vert T_N(E,N,0)\Vert$ is bounded, the result follows.

Now, by Proposition 3.9 in \cite{Lasttransfermatrix}, one has for every $E \in [0,1]$ and every $L\in\mathbb{N}$,
\begin{eqnarray*} \mathfrak{Im} \biggr(\int_{\mathbb{R}} \frac{1}{E+i/L - z} d\mu_{\delta_1}^{N}(z) \biggr) \leq (5 + \sqrt{24}) \frac{1}{L} \sum_{n=0}^{L+1} \|T_N(E,n,0)\|^2.
\end{eqnarray*}
Thus, for every $E \in [0,1]$ and every $L\ge N$,
\begin{eqnarray*}  \mathfrak{Im} \biggr(\int_{\mathbb{R}} \frac{1}{E+i/L - z} d\mu_{\delta_1}^{N}(z) \biggr) &\leq& (5 + \sqrt{24}) \frac{1}{L} \sum_{n=0}^{L+1} \|T_N(E,n,0)\|^2\\ &=& (5 + \sqrt{24}) \biggr( \frac{1}{L} \sum_{n=0}^{N} \|T_N(E,n,0)\|^2 + \frac{1}{L} \sum_{n=N+1}^{L+1} \|T_N(E,n,0)\|^2 \bigg)\\ &\leq& (5 + \sqrt{24}) \biggr( \frac{1}{L} \sum_{n=0}^{N} \|T_N(E,n,0)\|^2 + \frac{L-N}{L} C_N \bigg).
\end{eqnarray*}
By Stone's Formula, it follows that for each $E \in [0,1]$ 
\begin{eqnarray*}
\frac{d\mu_{\delta_1}^{N}}{dx}(E) &=& \frac{1}{\pi} \lim_{ L \to \infty} \mathfrak{Im} \biggr\langle \delta_1, R^N\biggr(E+i\frac{1}{L}\biggr)\delta_1 \biggr\rangle = \frac{1}{\pi}  \lim_{ L \to \infty} \mathfrak{Im} \biggr(\int_{\mathbb{R}} \frac{1}{E+i/L - z} d\mu_{\delta_1}^{N}(z) \biggr).
\end{eqnarray*}
Since $\sum_{n=0}^{N} \|T_N(E,n,0)\|^2$ does not depend on $L$, one gets 
\[C_{2,N}:= \frac{C_N}{\pi}(5+\sqrt{24}) \geq \displaystyle\sup_{E \in [0,1]} \frac{d\mu_{\delta_1}^{N}}{dx}(E).\]
\end{proof}

\subsection{Power-law singularities and quantum dynamics}

Recall that, for every $N \in \mathbb{N} \cup \{0\}$,  ${\mathbb{R}} \ni t \mapsto e^{-itH_N}$ is a one-parameter strongly continuous unitary evolution group and, for each $\psi\in \ell^2({\mathbb{N}})$, $(e^{-itH_N}\psi)_{t \in \mathbb{R}}$ is the unique solution to the Schr\"odinger equation
\[
\begin{cases} \partial_t \psi = -iH_N\psi, \quad t \in {\mathbb{R}}, \\ \psi(0) = \psi\in\ell^2(\mathbb{N}).  \end{cases}
\]

Now, we present a dynamical quantity usually considered to probe the large time behavior of $e^{-itH_N}\psi$, the  so-called (time-average) {\em quantum return probability}, which gives the (time-average) probability of finding the particle at time $t>0$ in its
initial state $\psi$:
\begin{equation*}
\frac{1}{t}\int_0^t |\langle \psi, e^{-isH_N} \psi \rangle|^2 \, ds.         
\end{equation*}
  
For  $0 \leq \beta < 1$ and  every $f \in {\mathrm L}^1(\mathbb{R})$, let
\[K_{\beta,f} = \{(|\cdot|^{-\beta}  \chi_{(0,1]}) \ast f\}.\]
Suppose that $f \geq 0$ and set 
\[\psi_{\beta,f} := \sqrt{K_{\beta,f}}(H_N)\delta_1,\]
where each $\sqrt{K_{\beta,f}}(H_N): \dom (\sqrt{K_{\beta,f}}(H_N)) \subset \ell^2({\mathbb{N}})\rightarrow \ell^2({\mathbb{N}})$ is defined through the functional calculus: for every $\psi\in \dom (\sqrt{K_{\beta,f}}(H_N))$, one has

\[\langle\psi,\sqrt{K_{\beta,f}}(H_N)\psi\rangle=\int \sqrt{K_{\beta,f}(x)}\,d\mu^{H_N}_\psi(x).\]

Note that $\delta_1 \in \dom (\sqrt{K_{\beta,f}}(H_N))$, since $\sqrt{K_{\beta,f}} \in {\mathrm L}^2(\mathbb{R},d\mu^{H_N}_{\delta_1})$; if  $\displaystyle\int_0^1 f(x) \, dx = 1$, then for every $x \in (0,1)$, $\displaystyle\lim_{\beta \downarrow 0} \sqrt{K_{\beta,f}}(x) = 1$, by dominated convergence; thus,  
\[\displaystyle\lim_{\beta \downarrow 0} \psi_{\beta,f} = \delta_1.\]
Our next result describes the behavior of the (time-average) {\em quantum return probability} of the initial states $\psi_{\beta,f}$. 

\begin{theorem}\label{maintheorem} Let $\psi_{\beta,f}$ be  as above, with $0\le f \in {\mathrm L}^1(\mathbb{R})$. Then:
\begin{enumerate}

\item[{\rm i)}] if $0\leq \beta < \frac{1}{2}$, for every $t>0$, 
\begin{eqnarray*}
\frac{1}{t}\int_0^t |\langle \psi_{\beta,f}, e^{-isH_N} \psi_{\beta,f} \rangle|^2 \, ds \leq   \frac{20\pi \|f\|_{{\mathrm L}^1(\mathbb{R})}^2}{(1-2\beta)}  \biggr\|\frac{d\mu_{\delta_1}^{N}}{dx}\biggr\|_{{\mathrm L}^\infty[0,1]}^2   t^{-1};
\end{eqnarray*}

\item[{\rm ii)}] if $\frac{1}{2}< \beta <1$, for every $t>0$, 
\begin{eqnarray*}
\frac{1}{t}\int_0^t |\langle \psi_{\beta,f}, e^{-isH_N} \psi_{\beta,f} \rangle|^2 \, ds   &\leq&  \frac{\Gamma(\beta - 1/2)  2^{18} e^{2\pi} {M_\beta} \|f\|_{{\mathrm L}^1(\mathbb{R})}^2  \biggr\|\frac{d\mu_{\delta_1}^{N}}{dx}\biggr\|_{{\mathrm L}^\infty[0,1]}^2 }{(2\pi)^{2(\beta-1)}} \;  t^{-2(1-\beta)};
\end{eqnarray*}

\item[{\rm iii)}] if $\beta = \frac{1}{2}$, for every $t>0$, 
\begin{eqnarray*}
\frac{1}{t}\int_0^t |\langle \psi_{\beta,f}, e^{-isH_N} \psi_{\beta,f} \rangle|^2 \, ds   &\leq& 4 \pi^2 e^{2\pi} \|f\|_{{\mathrm L}^1(\mathbb{R})}^2  \biggr\|\frac{d\mu_{\delta_1}^{N}}{dx}\biggr\|_{{\mathrm L}^\infty[0,1]}^2\biggr[ \frac{1}{t} +  M_{\frac{1}{2}}  \frac{\Gamma(0, 16\pi^4/t^2)}{t} \biggr],
\end{eqnarray*}
in particular, for sufficiently large $t$ 
\begin{eqnarray*}
\frac{1}{t}\int_0^t |\langle \psi_{\beta,f}, e^{-isH_N} \psi_{\beta,f} \rangle|^2 \, ds   &\leq& 4 \pi^2 e^{2\pi} \|f\|_{{\mathrm L}^1(\mathbb{R})}^2  \biggr\|\frac{d\mu_{\delta_1}^{N}}{dx}\biggr\|_{{\mathrm L}^\infty[0,1]}^2\biggr[ \frac{1}{t} +  3M_{\frac{1}{2}}  \frac{\log(t)}{t} \biggr].
\end{eqnarray*}
\end{enumerate}
\end{theorem}

\begin{remark}{\rm By applying Theorem \ref{2Stheorem}, in general, one can extend the above result to families of Schr\"odinger operators  whose Radon-Nikodym derivatives of each spectral measure (with respect to  Lebesgue measure) is bounded (see the proof of Theorem~\ref{maintheorem}).}
\end{remark}

We revisit Example \ref{ex1}, but now taking into account  Theorem~\ref{maintheorem}. 

\begin{example}\label{ex2} {\rm Let $\frac{1}{2} \leq \beta<1$ and   $0<\delta<1-\beta$. Let $f_\delta \in {\mathrm L}^1(\mathbb{R})$ with
 $\displaystyle\int_0^1 f_\delta(x) dx =1$ and suppose that there exists $C>0$ so that, for every $w \in (0,1]$, 
\[f_\delta(w) \geq \frac{C}{w^{1-\delta}}.\]
Hence, for every $0< w \leq 1$,
\begin{eqnarray*}
K_{\beta,f_\delta}(w) &\geq&  C\int_0^1 \frac{1}{y^\beta |w-y|^{1-\delta}} \, dy \geq  C\int_0^w \frac{1}{y^\beta |w-y|^{1-\delta}} \, dy \geq \frac{C}{w^\beta} \int_0^w \frac{1}{|w-y|^{1-\delta}} \, dy \geq   C w^{- \beta+\delta}.
\end{eqnarray*}
Thus, for every $0< \epsilon < 1$,
\[ \mu_{\psi_{\beta,f_\delta}}((0,\epsilon))=\int_0^\epsilon K_{\beta,f}(w) \frac{d\mu_{\delta_1}^{N}}{dx}(w) \, dw \geq C C_{1,N} \int_0^\epsilon w^{-\beta+\delta} dw = C C_{1,N} \frac{\epsilon^{(1-\beta+ \delta)} }{(1-\beta+ \delta)},\]
with $C_{1,N}$ given by Lemma \ref{teclemma}. Therefore, in this case, 
\[ d\mu_{\psi_{\beta,f_\delta}} = K_{\beta,f_\delta} \frac{d\mu_{\delta_1}^{N}}{dx} \, dx\]
is at most  {\rm U}$(1-{\beta}+ \delta)${\rm H}, and so, by Theorem~\ref{Strichartztheorem}-i), one can say at most that
\begin{eqnarray*}\frac{1}{t}\int_0^t |\langle \psi_{\beta,f_\delta}, e^{-isH_N} \psi_{\beta,f_\delta} \rangle|^2 \,ds = O(t^{-(1-\beta)-\delta}).
\end{eqnarray*}
Nonetheless, it follows from Theorem  \ref{maintheorem} that, for every $\frac{1}{2} < \beta < 1$,
\begin{equation}\label{eqex2}
  \frac{1}{t}\int_0^t \vert\langle \psi_{\beta,f_\delta}, e^{-isH_N} \psi_{\beta,f_\delta} \rangle\vert^2 \, ds = O(t^{-2(1-\beta)})
\end{equation}
and
\begin{equation}\label{eq9999}
  \frac{1}{t}\int_0^t \vert\langle \psi_{1/2,f_\delta}, e^{-isH_N} \psi_{1/2,f_\delta} \rangle\vert^2 \, ds = O(\log(t)/t))
\end{equation} 

We observe that the above rates are power-law optimal. Namely, if there exists $\varepsilon>0$ so that one can replace $t^{-2(1-\beta)}$ by $t^{-2(1-\beta)+\varepsilon}$  in~\eqref{eqex2} or (\ref{eq9999}), then by Theorem \ref{Strichartztheorem}-ii),  $\psi_{\beta,f_\delta}$ will be at least uniformly $(1-\beta+\varepsilon/4)$-H\"older; for $0<\delta< \varepsilon/4$, one gets a contradiction.} 
\end{example}

\begin{remark}{\rm For each $0<\delta<1-\beta$, Example \ref{ex2} presents a family of initial states, $\psi_{\beta,f_\delta}$, such that  $\displaystyle\lim_{\beta \downarrow 0} \psi_{\beta,f_\delta} = \delta_1$  and for which the correspondent (time-average) quantum return probabilities depend continuously on the power-law growth rates of the singularities of the respective spectral measures.}
\end{remark}

\begin{proof}[{Proof} {\rm (Theorem~\ref{maintheorem})}]  i) It follows from the  Spectral Theorem that for every $t>0$, 
\begin{eqnarray*}
\frac{1}{t}\int_0^t |\langle \psi_\beta, e^{-isH_N} \psi_\beta \rangle|^2 \, ds &=&  \frac{1}{t} \int_0^t \bigg|\int_{\mathbb{R}} e^{-isy} K_{\beta,f}(y) d\mu_{\delta_1}^{N}(y) \bigg|^2  ds \\ &=&  \frac{2\pi}{t} \int_0^{\frac{t}{2\pi}} \bigg|\int_{\mathbb{R}} e^{- 2\pi isy} K_{\beta,f}( y) \frac{d\mu_{\delta_1}^{N}}{dx}(y)dy \bigg|^2  ds.
\end{eqnarray*}
Since $0\leq \beta < \frac{1}{2}$, by Young's Convolution Inequality and Lemma \ref{teclemma}, $K_{\beta,f} \frac{d\mu_{\delta_1}^{N}}{dx} \in {\mathrm L}^2(\mathbb{R})$. Thus, by Theorem \ref{Strichartztheorem}-i) applied to $d\mu=\chi_{[0,1]} dy$ and to the function $K_{\beta,f}( y) \frac{d\mu_{\delta_1}^{N}}{dx}( y)$, one obtains, for every $t>0$, 
\begin{eqnarray*}
  \frac{1}{t}\int_0^t |\langle \psi_{\beta,f}, e^{-itH_N} \psi_{\beta,f} \rangle|^2 \, ds &=& \frac{2\pi}{t} \int_0^{\frac{t}{2\pi}} \bigg|\int_{\mathbb{R}} e^{- 2\pi isy} K_{\beta,f}( y) \frac{d\mu_{\delta_1}^{N}}{dx}(y)dy \bigg|^2  ds \\ &\leq&  10 \|K_{\beta,f}\|_{{\mathrm L}^2(\mathbb{R})}^2 \biggr\|\frac{d\mu_{\delta_1}^{N}}{dx}\biggr\|_{{\mathrm L}^\infty[0,1]}^2  2\pi t^{-1} \\
  &\leq &\frac{20\pi \|f\|_{{\mathrm L}^1(\mathbb{R})}^2}{(1-2\beta)}  \biggr\|\frac{d\mu_{\delta_1}^{N}}{dx}\biggr\|_{{\mathrm L}^\infty[0,1]}^2   t^{-1}.
\end{eqnarray*}
We remark that if $d\mu = \chi_{[0,1]} dx$, then one can choose $C_\mu =10$ in Theorem \ref{Strichartztheorem}-i); for details, see page~416 in \cite{Last}.

\noindent ii) and iii) Let $\frac{1}{2}\leq \beta < 1$; since for every $t>0$,

\begin{eqnarray*}
\frac{1}{t}\int_0^t |\langle \psi_{\beta,f}, e^{-isH_N} \psi_{\beta,f} \rangle|^2 \, ds  =   \frac{2\pi}{t} \int_0^{\frac{t}{2\pi}}  \bigg|\int_{\mathbb{R}} e^{-2\pi isy} d\mu_{\beta,f,g}(y) \bigg|^2  ds,  
\end{eqnarray*}
where $g = \displaystyle\frac{d\mu_{\delta_1}^{N}}{dx}$, ii) and iii) are direct consequences of Lemma \ref{teclemma} and Theorem \ref{2Stheorem}. 
\end{proof}


\section{Proof of Theorem \ref{2Stheorem}}\label{sectProofMain}

\begin{lemma}\label{mainlemma} Let $ K_{\beta,f}$ be as before.

\begin{enumerate}
\item[{\rm i)}] If $\frac{1}{2} < \beta < 1$, then, for every $t>0$, one has

\begin{eqnarray*}
\bigg| \int_{\mathbb{R}} \int_{\mathbb{R}} e^{-\pi t^2|x-y|^2} K_{\beta,f}(x) \overline{K_{\beta,f}} (y)  \, dx dy \bigg| \; \leq \; {M_\beta} \Gamma(\beta - 1/2) \|f\|_{{\mathrm L}^1(\mathbb{R})}^2\; t^{-2(1-\beta)}.
\end{eqnarray*}

\item[{\rm ii)}]  If $\beta = \frac{1}{2}$, then, for every $t>0$, one has

\begin{eqnarray*}
\bigg| \int_{\mathbb{R}} \int_{\mathbb{R}} e^{-\pi t^2|x-y|^2} K_{\beta,f}(x) \overline{K_{\beta,f}} (y)  \, dx dy \bigg| \; \leq \; \frac{\pi\|f\|_{{\mathrm L^1}(\mathbb{R})}^2}{t} +  M_{\frac{1}{2}} \|f\|_{{\mathrm L^1}(\mathbb{R})}^2 \frac{\Gamma(0, \pi/t^2)}{t},  
\end{eqnarray*}
and recall that $M_\beta$ is given by~\eqref{eqMbeta}.
\end{enumerate}
\end{lemma}

\begin{proof} Let $\frac{1}{2} \leq \beta <1$ and let $(K_n)_{n \in \mathbb{N}} \subset {\mathrm L}^1(\mathbb{R}) \cap {\mathrm L}^2(\mathbb{R})$ be so that $\displaystyle\lim_{n \to \infty}\|K_n - K_{\beta,f}\|_{{\mathrm L}^1(\mathbb{R})} = 0$. Then, by Theorem 4.9 in \cite{Brezis}, there exist a subsequence $(K_{n_k})$ and a function $h \in {\mathrm L}^1(\mathbb{R})$ such that $\displaystyle\lim_{k \to \infty} K_{n_k}(x) = K_{\beta,f}(x)$ for almost every $x \in \mathbb{R}$, and for every $k \geq 1$, $ |K_{n_k}(x)| \leq h(x)$ for  almost every $x \in \mathbb{R}$. We note that for each $t >0$, each $k \geq 1$ and each $\xi \in \mathbb{R}$, 

\[  e^{- \frac{\pi |\xi|^2}{t^2}} |\widehat{K_{n_k}}(\xi)|  \leq  e^{- \frac{\pi |\xi|^2}{t^2}} \|K_{n_k}\|_{{\mathrm L}^1(\mathbb{R})}  \leq  e^{- \frac{\pi |\xi|^2}{t^2}} \|h\|_{{\mathrm L}^1(\mathbb{R})},\]
and for each $t >0$, 
\[\int_{\mathbb{R}} e^{- \frac{\pi |\xi|^2}{t^2}} d\xi = t.\]
This show that, for every $t>0$, the sequence $e^{- \frac{\pi |\xi|^2}{t^2}} |\widehat{K_{n_k}}(\xi)|$ is dominated by a integrable function.  
 
Set, for each $t>0$ and each $x \in \mathbb{R}$, $\Phi_t(x) := e^{-\pi t|x|^2}$. Then, for each $t>0$,
\begin{equation}\label{eq0lemma}
\widehat{\Phi_t}(\xi) = \frac{1}{t} e^{- \frac{\pi |\xi|^2}{t^2}}, \quad \xi \in \mathbb{R}.
\end{equation}

It follows from the identity in (\ref{eq0lemma}), some basic properties of the Fourier transform, dominated convergence and Plancherel's Theorem that for each $y \in \mathbb{R}$ and each $t >0$, 
\begin{eqnarray}\label{eq00lemma}\nonumber  \int_{\mathbb{R}} e^{-\pi t^2|x-y|^2} K_{\beta,f}(x) \, dx  &=& \lim_{k \to \infty} \int_{\mathbb{R}} e^{-\pi t^2|x-y|^2} K_{n_k} (x)  \, dx \nonumber =  \lim_{k \to \infty} \int_{\mathbb{R}} \overline{(\tau_y \Phi_t)(x)} K_{n_k} (x)\, dx \\ \nonumber  &=& \lim_{k \to \infty} \int_{\mathbb{R}}  \overline{\widehat{(\tau_y \Phi_t)}(\xi)} \widehat{K_{n_k}}(\xi) \,  d\xi = \lim_{k \to \infty}  \int_{\mathbb{R}} e^{2\pi i y \xi} \widehat{\Phi_t}(\xi)      \widehat{K_{n_k}}(\xi) \, d\xi  \\    \nonumber &=& \lim_{k \to \infty} \frac{1}{t} \int_{\mathbb{R}} e^{2\pi i y \xi} e^{- \frac{\pi |\xi|^2}{t^2}} \widehat{K_{n_k}}(\xi)  \, d\xi\\ &=& \frac{1}{t} \int_{\mathbb{R}}  e^{2\pi i y \xi} e^{- \frac{\pi |\xi|^2}{t^2}}  \widehat{K_{\beta,f}}(\xi) \, d\xi,
\end{eqnarray}
where $\tau_yf(\cdot) = f(\cdot-y)$ stands for the translation by $y\in\mathbb{R}$. Thus, by Fubini's Theorem, one obtains for each $t >0$,   
\begin{eqnarray}\label{eq1lemma} \nonumber  \int_{\mathbb{R}} \int_{\mathbb{R}} e^{-\pi t^2|x-y|^2} K_{\beta,f}(x) \overline{K_{\beta,f}} (y)  \, dx dy   &=& \frac{1}{t} \int_{\mathbb{R}}  \int_{\mathbb{R}} e^{2\pi i y \xi}  \overline{K_{\beta,f}} (y) \, dy \, e^{- \frac{\pi |\xi|^2}{t^2}}  \widehat{K_{\beta,f}}(\xi) \,  d\xi\\ &=&  \frac{1}{t} \int_{\mathbb{R}}  e^{- \frac{\pi |\xi|^2}{t^2}}  |\widehat{K_{\beta,f}}(\xi)|^2   \, d\xi.  
\end{eqnarray}

Now, by the Convolution Theorem, it follows that for each $\xi >0$,  
\begin{eqnarray}\label{eq88888} \nonumber \widehat{K_{\beta,f}}(\xi) &=& \biggl\{\int_0^1 e^{-2\pi i x \xi } x^{-\beta} dx \biggl\} \, \hat{f}(\xi) = \biggl\{\frac{1}{(2\pi)^{1-\beta}\xi^{1-\beta}} \int_0^1 e^{-2\pi i x \xi } (2 \pi x \xi)^{-\beta} (2\pi \xi)dx \biggl\} \, \hat{f}(\xi)\\ &=& \biggl\{\frac{1}{(2\pi)^{1-\beta}\xi^{1-\beta}} \int_0^{2\pi\xi} e^{- iu} u^{-\beta} du \biggl\} \, \hat{f}(\xi).
\end{eqnarray}
One also has, for each $\xi <0$,   that 
\begin{eqnarray*} \widehat{K_{\beta,f}}(\xi) = \biggl\{\frac{1}{(2\pi)^{1-\beta}(-\xi)^{1-\beta}} \overline{\int_0^{-2\pi\xi} e^{-iu} u^{-\beta} du} \biggl\} \, \overline{\hat{f}(\xi)}.
\end{eqnarray*}
Thus, for  $\xi \neq 0$,
\begin{eqnarray}\label{eq2lemma} |\widehat{K_{\beta,f}}(\xi)|^2 \leq  {M_\beta}  \, \frac{\|f\|^2_{{\mathrm L^1}(\mathbb{R})}}{|\xi|^{2(1-\beta)}}.
\end{eqnarray}

Now a separate argument is  necessary for each item. i) $1/2<\beta<1$. Since it follows from Cauchy's Residue Theorem that for every $t >0$  
\begin{equation}\label{eq01}
\int_{\mathbb{R}} \frac{e^{- \frac{\pi |\xi|^2}{t^2}} }{|\xi|^{2(1-\beta)}  }  \, d\xi =  \pi^{1/2-\beta}\Gamma\biggr(\beta - \frac{1}{2}\biggr) t^{2\beta-1} \leq \Gamma\biggr(\beta - \frac{1}{2}\biggr) t^{2\beta-1} ,
\end{equation}
one gets from (\ref{eq1lemma}) and (\ref{eq2lemma}) that for every $t >0$,  
\begin{eqnarray*}
 \int_{\mathbb{R}} \int_{\mathbb{R}} e^{-\pi t^2|x-y|^2} K_{\beta,f}(x) \overline{K_{\beta,f}} (y)  \, dx dy     &=&  \frac{1}{t} \int_{\mathbb{R}}  e^{- \frac{\pi |\xi|^2}{t^2}}  |\widehat{K_{\beta,f}}(\xi)|^2   \, d\xi  \\  &\leq&   \frac{{M_\beta}  \|f\|_{{\mathrm L^1}(\mathbb{R})}^2}{t} \int_{\mathbb{R}} \frac{e^{- \frac{\pi |\xi|^2}{t^2}} }{|\xi|^{2(1-\beta)}  }  \, d\xi\\  &\leq& \Gamma(\beta - 1/2) {M_\beta}   \|f\|_{{\mathrm L^1}(\mathbb{R})}^2 t^{-2(1-\beta)}.
\end{eqnarray*}

ii) $\beta=1/2$. By (\ref{eq88888}), for every $t> 0$
\begin{eqnarray*} \int_{-1}^1  e^{- \frac{\pi |\xi|^2}{t^2}}  |\widehat{K_{\beta,f}}(\xi)|^2   \, d\xi  &\leq&  2\|f\|_{{\mathrm L^1}(\mathbb{R})}^2 \int_0^1  e^{- \frac{\pi |\xi|^2}{t^2}}  \frac{1}{\xi} \biggr(\int_0^{2\pi\xi} \frac{1}{\sqrt{u}} \, du \biggr)^2 \, d\xi \\ &=&   \pi \|f\|_{{\mathrm L^1}(\mathbb{R})}^2 \int_{0}^1  e^{- \frac{\pi |\xi|^2}{t^2}} \,  d\xi \leq  \pi\|f\|_{{\mathrm L^1}(\mathbb{R})}^2.
\end{eqnarray*}
Since, by Cauchy's Residue Theorem, for every $t >0$
\[\int_1^\infty \frac{e^{- \frac{\pi |\xi|^2}{t^2}} }{\xi}  \, d\xi = \frac{\Gamma(0, \pi/t^2)}{2},\]
one gets from (\ref{eq1lemma}) and (\ref{eq2lemma}) that for every $t >0$,
\begin{eqnarray*}
 \int_{\mathbb{R}} \int_{\mathbb{R}} e^{-\pi t^2|x-y|^2} K_{\beta,f}(x) \overline{K_{\beta,f}} (y)  \, dx dy     &=&  \frac{1}{t} \int_{\mathbb{R}}  e^{- \frac{\pi |\xi|^2}{t^2}}  |\widehat{K_{\beta,f}}(\xi)|^2   \,  d\xi  \\  &=& \frac{1}{t} \int_{-1}^1  e^{- \frac{\pi |\xi|^2}{t^2}}  |\widehat{K_{\beta,f}}(\xi)|^2   \, d\xi  + \frac{1}{t} \int_{|x|> 1}  e^{- \frac{\pi |\xi|^2}{t^2}}  |\widehat{K_{\beta,f}}(\xi)|^2   \, d\xi\\ &\leq&  \frac{1}{t} \int_{-1}^1  e^{- \frac{\pi |\xi|^2}{t^2}}  |\widehat{K_{\beta,f}}(\xi)|^2   \, d\xi +  \frac{ 2M_{\frac{1}{2}} \|f\|_{{\mathrm L^1}(\mathbb{R})}^2}{t} \int_1^\infty \frac{e^{- \frac{\pi |\xi|^2}{t^2}} }{\xi}  \, d\xi\\ &\leq&   \frac{\pi\|f\|_{{\mathrm L^1}(\mathbb{R})}^2}{t} +  M_{\frac{1}{2}} \|f\|_{{\mathrm L^1}(\mathbb{R})}^2 \frac{\Gamma(0, \pi/t^2)}{t}.
\end{eqnarray*} 
\end{proof}

\begin{remark}{\rm  Since the  integral $\displaystyle\int_{\mathbb{R}} \frac{e^{- \frac{\pi |\xi|^2}{t^2}} }{|\xi|  }  \, d\xi$  does not converge,   a separate argument is necessary therein for the case $\beta = \frac{1}{2}$.}
\end{remark}

\begin{proof} [{Proof} {\rm (Theorem~\ref{2Stheorem})}] Note that the proof of Theorem~\ref{2Stheorem} is a consequence of Lemma~\ref{mainlemma} and Fubini's Theorem. We present  details of the proof of  item~i): by the linearity of the convolution product, we may assume without loss of generality that $\|f\|_{{\mathrm L^1}(\mathbb{R})} \leq 1 $ and $ \|g\|_{{\mathrm L}^\infty[0,1]} \leq 1$.  We divide this proof into two cases.

\

\noindent {\bf Case 1:} $f,g$ are nonnegative real-valued functions. By Fubini's Theorem, one has, for every $t>0$,   

\begin{eqnarray}\label{maineq1}
\nonumber \frac{1}{t} \int_0^t \bigg|\int_{\mathbb{R}} e^{-2\pi isx} K_{\beta,f}(x) g(x) \, dx \bigg|^2  ds &\leq& \frac{1}{t} \int_0^t \bigg|\int_{\mathbb{R}} e^{-2\pi isx} K_{\beta,f}(x) g(x) \, dx \bigg|^2 \, e^{2\pi-(2\pi s)^2/t^2} \, ds\\ \nonumber &\leq& \frac{e^{2\pi}}{t} \int_{-\infty}^{\infty} \bigg|\int_{\mathbb{R}} e^{-2\pi isx} K_{\beta,f}(x) g(x) \, dx \bigg|^2 \, e^{-(2\pi s)^2/t^2} \, ds\\ \nonumber &=& \frac{e^{2\pi}}{t}  \int_{\mathbb{R}} \int_{\mathbb{R}} K_{\beta,f}(x) g(x) \overline{K_{\beta,f}(y)} \, \overline{g(y)}\\ \nonumber &\times&  \biggl\{ \int_{-\infty}^{\infty} \,  e^{-((2\pi s)^2/t^2)-2\pi is(x-y)}  \, ds \biggl\} dxdy\\ \nonumber &=&  \frac{e^{2\pi} \sqrt{\pi}}{2 \pi}  \int_{\mathbb{R}} \int_{\mathbb{R}} K_{\beta,f}(x) g(x)\, K_{\beta,f}(y) g(y)  e^{-\frac{t^2|x-y|^2}{4}} dx dy \\ \nonumber  &\leq&  \frac{e^{2\pi} \sqrt{\pi}}{2 \pi}  \int_{\mathbb{R}} \int_{\mathbb{R}} K_{\beta,f}(x) K_{\beta,f}(y)   e^{-\frac{t^2|x-y|^2}{4}} dx dy \\  &=&  \frac{e^{2\pi} \sqrt{\pi}}{2 \pi}  \int_{\mathbb{R}} \int_{\mathbb{R}} K_{\beta,f}(x) K_{\beta,f}(y)   e^{-\pi (t/2\sqrt{\pi})^2|x-y|^2} dx dy .
\end{eqnarray}

It then follows from (\ref{maineq1})  combined with Lemma \ref{mainlemma} i) that, for every $t>0$,
\begin{eqnarray*}
 \frac{1}{t} \int_0^t \bigg|\int_{\mathbb{R}} e^{-2\pi isx} K_{\beta,f}(x) g(x) \, dx \bigg|^2  ds &\leq&  \frac{\Gamma(\beta - 1/2)  e^{2\pi} {M_\beta} \sqrt{\pi}}{2^{2\beta-1} \pi^\beta}\;   t^{-2(1-\beta)}.
\end{eqnarray*}

\noindent {\bf Case 2:} $f,g$ are complex valued. This case is a  direct consequence of {\bf Case 1}. Namely, by the linearity of the convolution product, by the inequality $(a+b)^2 \leq 2(a^2+b^2)$, $a,b>0,$ and by the identity 
\begin{eqnarray*}
K_{\beta,f} \cdot g&=& \biggr\{K_{\beta,\mathfrak{Re}(f)^+}  - K_{\beta,\mathfrak{Re}(f)^-}  + i\biggr(K_{\beta,\mathfrak{Im}(f)^+} - K_{\beta,\mathfrak{Im}(f)^-}\biggr)\biggr\}\\ &\times& \biggr\{\mathfrak{Re}(g)^+ - \mathfrak{Re}(g)^- + i\biggr(\mathfrak{Im}(g)^+ - \mathfrak{Im}(g)^-\biggr)\biggr\},
\end{eqnarray*}
it follows that, for every $t>0$,
\begin{eqnarray*}
 \frac{1}{t} \int_0^t \bigg|\int_{\mathbb{R}} e^{-2\pi isx} K_{\beta,f}(x) g(x) \, dx \bigg|^2  ds &\leq&  \frac{\Gamma(\beta - 1/2) 2^{16}e^{2\pi} {M_\beta} \sqrt{\pi}}{2^{2\beta-1} \pi^\beta} \;  t^{-2(1-\beta)}.
\end{eqnarray*}
\end{proof}

 
\begin{center} \Large{Acknowledgments} 
\end{center}
\addcontentsline{toc}{section}{Acknowledgments}

\noindent M. Aloisio thank the partial support by CAPES (a Brazilian government agency; Finance Code 001). S. L. Carvalho thanks the partial support by FAPEMIG (Minas Gerais state agency; under contract 001/17/CEX-APQ-00352-17) and C. R. de Oliveira thanks the partial support by CNPq (a Brazilian government agency, under contract 303689/2021-8).


\noindent  Email: moacir@ufam.edu.br, Departamento de Matem\'atica, UFAM, Manaus, AM, 369067-005 Brazil

\noindent  Email: silas@mat.ufmg.br, Departamento de Matem\'atica, UFMG, Belo Horizonte, MG, 30161-970 Brazil

\noindent  Email: oliveira@ufscar.br,  Departamento  de  Matem\'atica,   UFSCar, S\~ao Carlos, SP, 13560-970 Brazil

\noindent  Email: edsonilustre@yahoo.com.br, Departamento de Matem\'atica, UFAM \& UEA, Manaus, AM, 369067-005 Brazil


\begin{thebibliography}{15}
\addcontentsline{toc}{section}{References}

\bibitem{AvilaUaH} A. Avila and S. Jitomirskaya, H\"older continuity of absolutely continuous spectral measures for one-frequency Schr\"odinger operators. Commun. Math. Phys. {\bf 301} (2011), 563--581.

\bibitem{Brezis} H. Brezis, Functional Analysis, Sobolev Spaces and PDE. Springer, (2011).

\bibitem{Damanik} D. Damanik and J. Fillman, Limit-periodic Schr\"odinger operators with Lipschitz continuous IDS. Proc. Am. Math. Soc. {\bf 147} (2019), 1531–1539.

\bibitem{Oliveira} C. R. de Oliveira, Intermediate Spectral Theory and Quantum Dynamics. Progress in Math. Phys. Basel, Birkh\"auser, (2009).

\bibitem{Germinet} F. Germinet, A. Kiselev and S. Tcheremchantsev, Transfer matrices and transport for
Schr\"odinger operators. Ann. Inst. Fourier. {\bf 54} (2004), 787–830

\bibitem{Last} Y. Last, Quantum dynamics and decompositions of singular continuous spectra. J. Funct. Anal. {\bf 42} (1996), 406--445

\bibitem{Lasttransfermatrix} Y. Last and B. Simon, Eigenfunctions, transfer matrices, and absolutely continuous spectrum of one-dimensional Schr\"odinger operators. Invent. Math. {\bf 135} (1999), 329–367.

\bibitem{Muller} V. M\"uller and Y. Tomilov, ``Large'' weak orbits of $C_0$-semigroups. Acta Sci. Math. (Szeged) {\bf 79} (2013), 475--505.

\bibitem{strichartz1990} R. S. Strichartz, Fourier asymptotics of fractal measures. J. Funct. Anal. {\bf89} (1990), 154--187.
 
\bibitem{Zhao} X. Zhao, H{\"o}lder continuity of absolutely continuous spectral measure for multi-frequency Schr{\"o}dinger operators. J. Funct. Anal. {\bf 278} (2020), 108508. 

\bibitem{Zhao2} X. Zhao, H{\"o}lder continuity of absolutely continuous spectral measure for the extended Harper's model. Nonlinearity. {\bf 34} (2021), 3356. 
 
\end{thebibliography}
\end{document}